\documentclass[10pt,a4paper]{article}    
\usepackage[margin=0.69in]{geometry}
\usepackage{authblk}
\usepackage{cite}
\usepackage{url}
\usepackage{tipa}
\usepackage[dvips]{graphicx}
\usepackage{amssymb}
\usepackage{color}
\usepackage{amsbsy}
\usepackage{textcomp}  
\usepackage{booktabs}
\usepackage{multirow}
\usepackage{marvosym}

\usepackage[cmex10]{amsmath}
\usepackage{graphicx}
\usepackage{amssymb}
\usepackage{mathrsfs}
\usepackage{cite}
\usepackage{bigints}

\newtheorem{lemma}{Lemma}
\newtheorem{remark}{Remark}

\newcommand{\bs}[1]{\boldsymbol{#1}}

\newcommand{\m}[1]{\mathcal{#1}}

\newcommand{\eor}{\hfill $\bigtriangledown$ \\} 

\newenvironment{myindentpar}[1]%
	{\begin{list}{}%
	{\setlength{\leftmargin}{#1}}%
	\item[]%
    			}
	{\end{list}}

\newcommand{\qed}{\nobreak \ifvmode \relax \else
      \ifdim\lastskip<1.5em \hskip-\lastskip
      \hskip1.5em plus0em minus0.5em \fi \nobreak
      \vrule height0.75em width0.5em depth0.25em\fi}
      
\usepackage{lipsum}

\title{\LARGE \textbf{A Data Augmentation Approach for a Class of Statistical Inference Problems}}

\author[$^\dagger$]{\Large Rodrigo Carvajal}
\author[$^\dagger$]{\Large Rafael Orellana}
\author[$^{\ast}$]{\Large Dimitrios Katselis}
\author[$^\ddagger$]{\Large Pedro Esc\'arate}
\author[$^{\dagger}$]{\Large Juan C. Ag\"uero}
\affil[$^\dagger$]{\normalsize Electronics Engineering Department, Universidad T\'{e}cnica Federico Santa Mar\'ia, Chile. \Letter~\{rodrigo.carvajalg,~juan.aguero\}@usm.cl, rafael.orellana.prato@gmail.com}
\affil[$^{\ast}$]{\normalsize Coordinated Science Laboratory, University of Illinois at Urbana-Champaign, IL, USA. \Letter~{katselis@illinois.edu}}
\affil[$^\ddagger$]{\normalsize Large Binocular Telescope Observatory, Steward Observatory, University of Arizona, Tucson, AZ, USA. \Letter~{pescarate@lbto.org}}
\date{}

\begin{document}
\maketitle

\begin{abstract}                        

We present an algorithm for a class of statistical inference problems. The main idea is to reformulate the inference problem as an optimization procedure, based on the generation of surrogate (auxiliary) functions. This approach is motivated by the MM algorithm, combined with the systematic and iterative structure of the Expectation-Maximization algorithm. The resulting algorithm can deal with hidden variables in Maximum Likelihood and Maximum a Posteriori estimation problems, Instrumental Variables, Regularized Optimization and Constrained Optimization problems.

The advantage of the proposed algorithm is to provide a systematic procedure to build  surrogate functions for a class of problems where hidden variables are usually involved. Numerical examples show the benefits of the proposed approach.  

\end{abstract}
\section{Introduction}
\label{sec:Intro}

Problems in statistics and system identification often involve variables for which measurements are not available. Among others, real-life examples can be found in communication systems \cite{Carvajal2013, Goldsmith2005} and systems with quantized data \cite{Godoy2011, Wang2003}. In Maximum Likelihood (ML) estimation problems, the \textit{likelihood function} is in general difficult to optimize by using closed-form expressions, and numerical approximations are usually cumbersome. These difficulties are traditionally avoided by the utilization of the Expectation-Maximization (EM) algorithm \cite{Dempster1977}, where a surrogate (auxiliary) function is optimized instead of the main objective function. This surrogate function includes the complete data, i.e. the measurements and the random variables for which there are no measurements. The incorporation of such \textit{hidden data} or \textit{latent variables} is usually termed as \textit{data augmentation}, where the main goal is to obtain, in general, simple and fast algorithms \cite{ref:Van_Dyk2001}.
 
On the other hand, the MM\footnote{MM stands for Maximization-Minorization or Minimization-Majorization, depending on the optimization problem that needs to be solved.} algorithm \cite{Hunter2004} is generally employed for solving more general optimization problems, not only for ML and Maximum a Posteriori (MAP) estimation problems. In general, the main motivation for using the MM algorithm is the lack of closed-form expressions for the solution of the optimization problem or dealing with objective cost functions that are not convex. Applications where the MM algorithm has been utilized include communication systems problems \cite{Pollakis2012} and image processing \cite{Figueiredo2007}. For constrained optimization problems, an elegant solution is presented by Marks and Wright \cite{Marks1978}, where the constraints are incorporated via the formulation of surrogate functions. Surprisingly, Marks' approach has not received the same attention from the scientific community when it comes to compare it with the EM and the MM algorithms. In fact, these three approaches are contemporary, but the EM algorithm has attracted most of the attention (out of the three methods), and it has been used for solving linear and nonlinear statistical inference problems in biology and engineering, see e.g. \cite{ref:Aguero2012, ref:Beal2003, ref:Gopaluni2008, ref:Hobolth2005, ref:Schon2010, ref:Yang2013}, amongst others. On the other hand, as shown in \cite{Hunter2004}, the MM algorithm has obtained much less attention, while Marks' approach is mostly known to a limited audience in the the Communication Systems community.  These three approaches have important similarities: i) a surrogate function is defined and optimized in place of the original optimization problem, and ii) the solution is obtained iteratively. In general, these algorithms are ``principles and recipes'' \cite{Meng2007} or a ``philosophy'' \cite{Hunter2004} for constructing solutions to a broad variety of optimization problems.

In this paper we adopt the ideas behind \cite{Dempster1977,Hunter2004,Marks1978} to develop an algorithm for a special class of functions. Our approach generalizes the ones of \cite{Dempster1977,Hunter2004,Marks1978} by reinterpreting the E-step in the EM algorithm and expressing the cost function in terms of an infinite mixture or kernel. In particular cases, the kernel corresponds to a variance-mean Gaussian mixture (VMGM), see e.g. \cite{Polson2013}. VMGMs, also referred to as normal variance-mean mixtures \cite{Barndorff-Nielsen1982} and normal scale mixtures \cite{West1987}, have been considered in the literature for formulating EM-based approaches to solve ML \cite{Balakrishnam2009} and MAP problems \cite{Polson2013, ref:Carvajal2015, Godoy2014}. Our approach is applicable to a wide class of functions, which allows for defining the likelihood function, the prior density function, and constraints as kernels, extending also the work in \cite{Marks1978}. Thus, our work encompasses the following contributions: i) a systematic approach to constructing surrogate functions for a class of cost functions and constraints, ii) a class of kernels where the unknown quantities of the algorithm can be easily computed, and iii) a generalization of \cite{Dempster1977, Hunter2004, Marks1978} by including the cost function and the constraints in one general expression. Our proposal is based, among other things, on a particular way to apply Jensen's inequality \cite[pp. 24--25]{Durret2010}. In addition, we provide the details on how to construct quadratic surrogate functions for cost functions and constraints.

Our algorithm is tested by two examples. In the first one we considered the problem of estimating the rotational velocities of stars. The system model corresponds to the convolution of two probability density functions (pdf's) and thus is an infinite mixture. We show that our reinterpretation of the EM algorithm allows for the direct application of our proposal for the correct estimation of the parameter of interest. In the second example we considered a system modelled by a linear regression. In this example the regression matrix is known and the parameter vector is unknown. The problem is solved considering a constrained optimization problem, where the (inequality) constraint is given by a multivariate $\ell_q$-norm of the parameter vector. In this example, we show that although the problem cannot be solved using the EM algorithm, the solution given by our proposal can iteratively converge to the global optimum.

\section{Rudiments of the proposed approach}
\label{section:bases}
\subsection{The EM algorithm}
\label{section:EM}

Let us consider an estimation problem and its corresponding log-likelihood function defined as $\ell(\bs \theta) = \log p(\bs y | \bs \theta)$, where $p(\bs y | \bs \theta)$ is the likelihood function, $\bs \theta \in \mathbb{R}^p$, and $\bs y \in \mathbb{R}^N$. Denoting the \textit{complete data} by $\bs z \in \Omega(\bs y)$, and using Bayes' theorem, we can obtain:
\begin{equation}
\ell(\boldsymbol{\theta}) = \log p(\bs y \vert \boldsymbol{\theta}) = \log p(\bs z\vert \boldsymbol{\theta}) - \log p(\bs z \vert \bs y,\boldsymbol{\theta}).
\label{eq:EM1}
\end{equation}
Let us assume that at the $i$th iteration we have the estimate $\hat{\boldsymbol{\theta}}^{(i)}$. By integrating at both sides of \eqref{eq:EM1} with respect to $p(\bs z\vert \bs y,\hat{\boldsymbol{\theta}}^{(i)})$ we obtain $ \ell(\boldsymbol{\theta})  
= \mathcal{Q}(\boldsymbol{\theta},\hat{\boldsymbol{\theta}}^{(i)} ) - \mathcal{H}(\boldsymbol{\theta},\hat{\boldsymbol{\theta}}^{(i)} )$,
where
\begin{align}
\mathcal{Q}(\boldsymbol{\theta},\hat{\boldsymbol{\theta}}^{(i)} ) = & \int_{\Omega(\bs y)} \log p(\bs z\vert\boldsymbol{\theta}) p(\bs z\vert \bs y,\hat{\boldsymbol{\theta}}^{(i)}) d\bs z, \label{eq:Q}\\
\mathcal{H}(\boldsymbol{\theta},\hat{\boldsymbol{\theta}}^{(i)} ) =&  \int_{\Omega(\bs y)} \log  p(\bs z\vert \bs y ,\boldsymbol{\theta})  p(\bs z\vert \bs y,\hat{\boldsymbol{\theta}}^{(i)}) d\bs z.
\label{eq:H}
\end{align}

Using Jensen's inequality \cite[pp. 24--25]{Durret2010}, it is possible to show that for any value of $\bs{\theta}$, the function  $\mathcal{H}(\boldsymbol{\theta},\hat{\boldsymbol{\theta}}^{(i)} )$ is  decreasing. Hence, the optimization is only carried out on the auxiliary function $\m{Q}({\bs{\theta}},\hat{{\bs{\theta}}}^{(i)})$ because, by maximizing $\m Q({\bs{\theta}},\hat{{\bs{\theta}}}^{(i)})$, the new parameter $\hat{\bs{\theta}}^{(i+1)}$ is such that the likelihood function increases (see e.g. \cite{Dempster1977,McLachlan1997}).

In general, the EM method can be summarised as follows:
\begin{myindentpar}{0cm}
\underline{E-step}: Compute the expected value of the joint \textit{likelihood function} 
for the \textit{complete data} (measurements and hidden variables) based on a given parameter estimate, $\hat{\bs{\theta}}^{(i)}$. Thus, we have (see e.g. \cite{Dempster1977}):
\begin{equation}
\mathcal{Q}({\bs{\theta}},\hat{{\bs{\theta}}}^{(i)}) = E[\,\log p(\bs z \vert {\bs{\theta}}) \vert \bs y,\hat{{\bs{\theta}}}^{(i)}\,],
\label{eq:E}
\end{equation}
\underline{M-step}: Maximize the function $\m{Q}({\bs{\theta}},\hat{{\bs{\theta}}}^{(i)})$ \eqref{eq:E}, with respect to $\bs{\theta}$:
\begin{equation}
\hat{\boldsymbol{\theta}}^{(i+1)} = \mathrm{arg} \max_{\boldsymbol{\theta}} \mathcal{Q}(\boldsymbol{\theta},\hat{\boldsymbol{\theta}}^{(i)} ).
\label{eq:M}
\end{equation}
\end{myindentpar}
This succession of estimates converges to a stationary point of the \textit{log-likelihood} function \cite{Vaida2005}.

\subsection{The MM algorithm}
\label{section:MM}

The idea behind the MM algorithm \cite{Hunter2004} is to construct a surrogate function $g({\bs{\theta}},\hat{{\bs{\theta}}}^{(i)})$, 
that majorizes (for minimization problems) or minorizes (for maximization problems) a given cost functions $f(\bs \theta)$ \cite{Hunter2004} at $\hat{{\bs{\theta}}}^{(i)}$ such that,
\begin{equation*}
\begin{array}{ll}
f(\bs \theta) \leq g(\bs \theta, \hat{\bs \theta}^{(i)}) & \qquad \text{for minimization problems, or} \\
f(\bs \theta) \geq g(\bs \theta, \hat{\bs \theta}^{(i)}) & \qquad \text{for maximization problems, and}\\
f(\bs \theta) = g(\bs{ \theta}^{(i)}, \bs{ \theta}^{(i)}), &
\end{array} 
\end{equation*}
where $\hat{\bs \theta}^{(i)}$ is an estimate of $\bs \theta$. Then, the surrogate function is iteratively optimized until convergence. Hence, for maximizing $f(\bs \theta)$ we have \cite{Sriperumbudur2011}
\begin{equation}
\hat{\bs\theta}^{(i+1)} = \text{arg} \, \max_{\bs \theta} g(\bs \theta, \hat{\bs\theta}^{(i)}).
\end{equation}
For the construction of the surrogate function, popular techniques include the second order Taylor approximation, the quadratic upper bound principle and Jensen's inequality for convex functions, see, e.g., \cite{Sriperumbudur2011}.
\begin{remark}
The iterative strategy utilized in the MM algorithm converges to a local optimum since 
\begin{equation*}
f(\hat{\bs \theta}^{(i+1)}) \geq g(\hat{\bs \theta}^{(i+1)},\hat{\bs \theta}^{(i)}) \geq g(\hat{\bs \theta}^{(i)},\hat{\bs \theta}^{(i)}) = f(\hat{\bs \theta}^{(i)}). 
\end{equation*} \eor
\end{remark}

\subsection{Data Augmentation in Inference Problems}

Data augmentation algorithms are based on the construction of the \textit{augmented data} and its many-to-one mapping $\Omega(y)$. This \textit{augmented data} is assumed to describe a model from which the observed data $y$ is obtained via marginalization \cite{ref:Tanner}. That is, a system with a likelihood function $p(y|\bs \theta)$ can be understood to arise from
\begin{equation}
p(y|\bs \theta) = \int p(y,x|\bs \theta)dx,
\end{equation}
where the \textit{augmented data} corresponds to $(y,x)$ and $x$ is the \textit{latent data} \cite{ref:Van_Dyk2001, ref:Tanner}.
This idea has been utilized for supervised learning \cite{ref:Figueiredo2003} and the development of the \textit{Bayesian Lasso} \cite{ref:Park2008}, to mention a few examples. In those problems, the Laplace distribution is expressed as a two-level hierarchical-Bayes model. This equivalence is obtained from the representation of the Laplace distribution as a VMGM:
\begin{equation}
\frac{a}{2} e^{-a|\theta|} = \bigintssss_0^\infty \underbrace{\frac{1}{\sqrt{2\pi \lambda}} e^{-\theta^2/(2\lambda)}}_{p(\theta|\lambda)}\underbrace{\frac{a^2}{2}e^{-a^2 \lambda/2}}_{p(\lambda)}d\lambda.
\end{equation}
In fact, there are several pdf's than can be expressed as VMGMs, as shown in Table \ref{table:VMGM} \cite{Polson2013}, where $g(\theta)$ is the penalty term that can be expressed as a pdf. In addition, in \cite{ref:Carvajal2015} it was developed an early version of the methodology presented in this paper, exploring the estimation of a  sparse parameter vector utilizing the $\ell_q$-(pseudo)norm, with $0<q<1$.
\begin{table}[h!]
\caption{Selection of mean-variance mixture representations for
  penalty functions. $p(\theta_j)\; =\; \int_0^{\infty} \m{N}_{\theta_j} (\mu_j \;+\; \lambda_j u_j,\tau^2 s_j^2 \lambda_j)p(\lambda_j)d\lambda_j$}   
\centering                          
\setlength{\tabcolsep}{2.5pt}
\begin{tabular}{c c c c c}             
\toprule
Penalty function & $g(\theta_j)$ & ${u}_j$ & $\mu_j$ & $p(\lambda_j)$  \\ [0.5ex]   
\cmidrule(r){1-5}                         
Ridge                     & $ (\theta_j/\tau)^2 $      								& $ 0 $ & $ 0 $ & $ \lambda_j = 1 $   \\
Lasso                     & $ |\theta_j/\tau|  $       							    & $ 0 $ & $ 0 $ & Exponential   \\
Bridge                    & $ |\theta_j/\tau|^\alpha $                              & $ 0 $ & $ 0 $ & Stable   \\
Generalized & \multicolumn{1}{c}{\multirow{2}{*}{$ \left [\frac{(1+\alpha)}{\tau} \right] \log \left ( 1 + \frac{|\theta_j|}{(\alpha\tau)}\right) $}}  & \multicolumn{1}{c}{\multirow{2}{*}{$ 0 $}} & \multicolumn{1}{c}{\multirow{2}{*}{$ 0 $}} & \multicolumn{1}{c}{\multirow{2}{*}{Exp-Gamma}}   \\
Double-Pareto & & & &  \\
\bottomrule                             
\end{tabular}
\label{table:VMGM}          
\end{table}

\section{Marks' approach for constrained optimization}
\label{section:Marks}

\subsection{Constrained problems in Statistical Inference}

Statistical Inference and System Identification techniques include a variety of methods that can be used in order to obtain a model of a system from data. Classical methods, such as \textit{Least Squares}, ML, MAP \cite{ref:Goodwin1977}, \textit{Prediction Error Method}, \textit{Instrumental Variables} \cite{ref:Soderstrom1988}, and \textit{Stochastic Embedding} \cite{ref:Ljung2014} have been considered in the literature for such task. However, the increasing complexity of modern system models has motivated researchers to revisit and reconsider those techniques for some problems. This has resulted in the incorporation of constraints and penalties, yielding an often more complicated optimization problem. For instance, it has been shown that the incorporation of linear equality constraints may improve the accuracy of the parameter estimates, see e.g \cite{ref:Mahata2004}. On the other hand, the incorporation of regularization terms (or penalties) also improves the accuracy of the estimates, reducing the effect of noise and eliminating spurious local minima \cite{ref:Lange2013}. Regularization can be mainly incorporated in two ways: by adding regularizing constraints (a penalty function) or by including a probability density function (pdf) as a prior distribution for the parameters, see e.g. \cite{Godoy2014}. Another way to improve the estimation is by incorporating inequality constraints, where certain functions of the parameters may be required, for physical reasons amongst others, to lie between certain bounds \cite{ref:Hanson1965}. From this point of view, it is possible to consider the classical methods with constraints or penalties, as in \cite{ref:Ljung2014, ref:Lange2013, ref:Hanson1965, ref:Hyder2009}.

Perhaps one of the most utilized approaches for penalized estimation (with complicated non-linear expressions) is the MM algorithm -- for details on the MM algorithm see Section \ref{section:MM}. This technique allows for the utilization of a surrogate function that is simple to handle, in terms of derivatives and optimization techniques, and that is, in turn, iteratively solved. However, its inequality constraint counterpart, here referred to as Marks' approach \cite{Marks1978}, is not as well known as the MM algorithm. Moreover, there is no straightforward manner to obtain such surrogate function. In this paper we focus on a systematic way to obtain the corresponding surrogate function using Marks' approach for a class of constraints. 

\subsection{Mark's approach}

The approach in \cite{Marks1978} deals with inequality constraints by using a similar approach to the EM and MM algorithms. The basic idea is, again, to generate a surrogate function that allows for an iterative procedure whose optimum value is the optimum value of the original optimization problem.

Let us consider the following constrained optimization problem:
\begin{align}
\bs \theta^{*} = & \, \text{arg} \min_{\bs \theta}  f(\bs \theta) \nonumber \\
\text{s. t.   }  & \,\, g(\bs{\theta}) \leq 0,
\label{eq:Opt_const}
\end{align}
where $f(\bs \theta)$ is the objective function and $g(\bs \theta)$ encodes the constraint of the optimization problem. In particular, let us focus on the case where $g(\bs \theta)$ is not a convex function. This implies that 
the optimization problem cannot be solved directly using standard techniques, such as quadratic programming or fractional programming. This difficulty can be overcome by utilizing a surrogate function $\m Q(\bs \theta, \hat{\bs \theta}^{(i)})$ at a given estimate $\hat{ \bs \theta} ^{(i)}$, such that
\begin{align}
g(\bs \theta) &\leq \m Q(\bs \theta, \hat{\bs \theta}^{(i)}) \label{eq:Marks_tight} \\
g(\hat{\bs \theta}^{(i)}) &= \m Q(\hat{\bs \theta}^{(i)}, \hat{\bs \theta}^{(i)}) \label{eq:Marks_decr} \\
\left. \frac{d }{d \bs{\theta}} g(\bs \theta) \right \vert_{\bs \theta = \hat{\bs \theta}^{(i)}} &= \left. \frac{d }{d \bs{\theta}} \m Q(\bs \theta, \hat{\bs \theta}^{(i)}) \right \vert_{\bs \theta = \hat{\bs \theta}^{(i)}} \label{eq:Marks_Fisher}
\end{align}
%
Provided the above properties are satisfied, then the following approximation of \eqref{eq:Opt_const}:
\begin{align}
\bs \theta^{(i+1)} = & \, \text{arg} \min_{\bs \theta}  f(\bs \theta) \nonumber \\
\text{s. t.   }  & \,\, \m Q(\bs{\theta},\hat{\bs{\theta}}^{(i)} ) \leq 0,
\label{eq:Marks}
\end{align}
iteratively converges to the solution of the optimization problem \eqref{eq:Opt_const}. As shown in \cite{Marks1978}, the optimization problem in \eqref{eq:Marks} is equivalent to the original problem \eqref{eq:Opt_const}, since the solution of \eqref{eq:Marks} converges to a point that satisfies the Karush-Kuhn-Tucker conditions of the original optimization problem.

\begin{remark}
Mark's approach can be considered as a generalization of the MM algorithm, since the latter can be derived (for a broad class of problems) from the former.
Let us consider the following problem:
\begin{equation}
\bs \theta ^\ast = \text{arg} \min_{\bs \theta}  f(\bs \theta).
\label{eq:unconst}
\end{equation}
Using the epigraph representation of \eqref{eq:unconst} \textup{\cite{ref:Grant2008}}, we obtain the equivalent problem
\begin{align}
\bs \theta ^\ast = &\,\,\text{arg} \min_{\bs \theta} \, t \nonumber \\
\text{s. t.   }  & \,\, f(\bs \theta) \leq t,
\end{align}
Using Mark's approach \eqref{eq:Marks}, we can iteratively find a local optimum of \eqref{eq:unconst} via
\begin{align}
\bs \theta^{(i+1)} = & \, \,\text{arg} \min_{\bs \theta} \, t \nonumber \\
\text{s. t.   }  & \,\, \m Q(\bs{\theta},\hat{\bs{\theta}}^{(i)} ) \leq t,
\label{eq:Marks_MM}
\end{align}
where $\m Q(\bs{\theta},\hat{\bs{\theta}}^{(i)} )$ in \eqref{eq:Marks_MM} is a surrogate function for $f(\bs \theta)$ in \eqref{eq:unconst}. From the epigraph representation we then obtain
\begin{equation}
\bs \theta^{(i+1)} =  \, \,\text{arg} \min_{\bs \theta} \, \m Q(\bs{\theta},\hat{\bs{\theta}}^{(i)} ),
\end{equation}
which is the definition of the MM algorithm (see \ref{section:MM}) for more details. \eor
\end{remark}

\section{A systematic approach to construct surrogate functions for a class of optimization problems}
\label{section:Gen_app}
Here, we consider a general optimization cost 
defined as:
\begin{equation}
\m{V}(\bs{\theta}) = \int \limits_{ \Omega(\bs{y})} K(\bs z, \bs \theta) d\mu(\bs z),
\label{eq:cost_kernel}
\end{equation}
where $\bs \theta$ is a parameter vector, $\bs y$ is a given data (i.e. measurements), $\bs z$ is the \textit{complete data} (comprised of the \textit{observed data} $\bs y$ and the \textit{hidden variables} (unobserved data), $\Omega(\bs{y})$ is a mapping from the sample space of $\bs z$ to the sample space of $\bs y$, $K(\cdot,\cdot)$ is a (positive) kernel function, and $\mu(\cdot)$ is a measure, see e.g \cite{Durret2010}. The definition in \eqref{eq:cost_kernel} is based on the definition of the auxiliary function $\m Q$ in the EM algorithm \cite[eq. (1.1)]{Dempster1977}, where it is assumed throughout the paper that there is a mapping that relates the \textit{not observed data} to the \textit{observed data}, and that the \textit{complete data} lies in $\Omega(\bs y)$ \cite{Dempster1977}. Notice that in \eqref{eq:cost_kernel} the kernel function may not be a pdf. However, several functions can be expressed in terms of a pdf. The most common cases are Gaussian kernels (yielding VMGMs) \cite{Barndorff-Nielsen1982} and Laplace kernels (yielding Laplace mixtures) \cite{ref:Garrigues2010}. 
\begin{remark}
Notice that, as explained in Section \ref{sec:Intro}, once the \textit{hidden data} has been selected, the \textit{data augmentation} procedure comes with the definition of $\m{V}(\bs{\theta})$ in \eqref{eq:cost_kernel}. From here, we follow the systematic nature of the EM and MM algorithms in terms of the iterative nature of the technique. \eor
\end{remark}

Since we are considering the optimization of the function $\m{V}(\bs{\theta})$, we can also consider the optimization of the function 
\begin{equation}
\m{J}(\bs \theta) = \log \m{V}(\bs{\theta}).
\label{eq:log_J}
\end{equation}
Without modifying the cost function in \eqref{eq:log_J}, we can multiply and divide by the logarithm of the kernel function, obtaining:
\begin{equation}
\m{J}(\bs \theta) = \log \m{V}(\bs{\theta}) = \log \m{V}(\bs{\theta}) \frac{\log K(\bs z,\bs \theta)}{\log K(\bs z,\bs \theta)} = \log K(\bs z,\bs \theta) - \log \frac{K(\bs z,\bs \theta)}{\m{V}(\bs{\theta})}.
\label{eq:log_J_prod}
\end{equation}
Let us assume that at the $i$th iteration we have the estimate ${\hat{\bs \theta}}^{(i)}$. Then, we can multiply by $\dfrac{K(\bs z,\hat{\bs \theta}^{(i)})}{\m{V}(\hat{\bs \theta}^{(i)}}$ and integrate on both sides of \eqref{eq:log_J_prod} with respect to $d\mu(\bs z)$, obtaining:
\begin{align}
\m{J}(\bs \theta) & = \int \limits_{ \Omega(\bs{y})}  \log \m{V}(\bs{\theta}) \frac{K(\bs z,\hat{\bs \theta}^{(i)})}{\m{V}(\hat{\bs \theta}^{(i)}} d\mu(\bs z) =  \log \m{V}(\bs{\theta}) \nonumber \\
&= \m{Q}(\bs \theta, \hat{\bs{\theta}}^{(i)}) - \m{H}(\bs \theta, \hat{\bs{\theta}}^{(i)}).
\label{eq:J_aux}
\end{align}
where :
\begin{align}
\m{Q}(\bs \theta, \hat{\bs{\theta}}^{(i)}) = & \int \limits_{ \Omega(\bs{y})} \log [K(\bs z, \bs \theta)] \frac{K(\bs z, \hat{\bs \theta}^{(i)})}{\m V(\hat{\bs \theta}^{(i)})} d\mu(\bs z), 
\label{eq:Q_gen}\\
\m{H}(\bs \theta, \hat{\bs{\theta}}^{(i)}) = & \int \limits_{ \Omega(\bs{y})} \log \left[ \frac{K(\bs z, \bs \theta)}{\m V({\bs \theta})}\right] \frac{K(\bs z, \hat{\bs \theta}^{(i)})}{\m V(\hat{\bs \theta}^{(i)})} d\mu(\bs z), 
\label{eq:H_gen}
\end{align}
are auxiliary functions.
%
As in the EM algorithm, for any $\bs \theta$, and using Jensen's inequality \cite[pp. 24--25]{Durret2010}, we have:
\begin{align}
\mathcal{H}(\boldsymbol{\theta},\hat{\boldsymbol{\theta}}^{(i)} ) \! - \mathcal{H}(\hat{\boldsymbol{\theta}}^{(i)},\hat{\boldsymbol{\theta}}^{(i)} ) \! & = \!\!\!\! \int \limits_{ \Omega(\bs{y})} \!\!\! \log \!\!\left[ \frac{K(\bs z, \bs \theta)}{\m V({\bs \theta})}\right] \!\! \frac{K(\bs z, \hat{\bs \theta}^{(i)})}{\m V(\hat{\bs \theta}^{(i)})} d\mu(\bs z) \nonumber\\
& - \int \limits_{ \Omega(\bs{y})} \log \left[ \frac{K(\bs z, \hat{\bs \theta}^{(i)})}{\m V(\hat{\bs \theta}^{(i)})}\right] \frac{K(\bs z, \hat{\bs \theta}^{(i)})}{\m V(\hat{\bs \theta}^{(i)})} d\mu(\bs z)\nonumber \\
& = \int \limits_{ \Omega(\bs{y})} \log \left[ \frac{K(\bs z, {\bs \theta})\m V(\hat{\bs \theta}^{(i)})}{\m V({\bs \theta})K(\bs z, \hat{\bs \theta}^{(i)})} \right] \frac{K(\bs z, \hat{\bs \theta}^{(i)})}{\m V(\hat{\bs \theta}^{(i)})} d\mu(\bs z)\nonumber \\
& \leq \log  \int \limits_{ \Omega(\bs{y})}  \frac{K(\bs z, \bs \theta)}{\m V({\bs \theta}}d\mu(\bs z) \nonumber\\
& = 0.
\label{eq:H_decreasing}
\end{align}
Hence, for any value of $\bs{\theta}$, the function  $\mathcal{H}(\boldsymbol{\theta},\hat{\boldsymbol{\theta}}^{(i)} )$ in \eqref{eq:H_gen}  is a decreasing function.

\begin{remark}
The kernel function $K(\bs z, \bs{\theta})$ satisfies the standing assumption $ K(\bs z, \bs{\theta}) > 0$ since the proposed scheme is built, among others, on the logarithm of the kernel function $K(\bs z, \bs{\theta})$. The definition of the kernel function in \eqref{eq:cost_kernel} allows for kernels that are not pdf's. On the other hand, some kernels may correspond to a scaled version of a pdf. In that sense, for the cost function in \eqref{eq:cost_kernel} we can define a new kernel and a new measure as
\[ \bar{K}(\bs z, \bs \theta)  = \frac{{K}(\bs z, \bs \theta)}{\int {K}(\bs z, \bs \theta) d\bs \theta}, \,\,\bar{\mu}(\bs z) = \left( \int {K}(\bs z, \bs \theta) d\bs \theta\right) {\mu}(\bs z),\]
\begin{equation*}
\Rightarrow \m{V}(\bs{\theta}) = \int \limits_{ \Omega(\bs{y})} \bar{K}(\bs z, \bs \theta) d\bar{\mu}(\bs z). 
\end{equation*} \eor
\end{remark}
\begin{remark}
In the proposed methodology, it is possible to optimize the surrogate function defined by
\begin{equation}
\bar{\m Q}(\bs \theta, \hat{\bs \theta }^{(i)} ) = \int \limits_{ \Omega(\bs{y})} \log [K(\bs z, \bs \theta)] {K(\bs z, \hat{\bs \theta}^{(i)})} d\mu(\bs z),
\label{eq:Qbar}
\end{equation}
since ${\m V(\hat{\bs \theta}^{(i)})}$ in \eqref{eq:Q_gen} does not depend on the parameter $\bs \theta$. 
Thus, the proposed method corresponds to a variation of the EM algorithm that is not limited to probability density functions (e.g. the likelihood function) for solving ML and MAP estimation problems. Instead, our version considers general measures ($\mu(\bs z)$), where the mapping over the measurement data $\Omega(\bs y)$ is a given set.\eor 
\end{remark}
The idea behind using a surrogate function is to obtain a simpler algorithm for the optimization of the objective function when compared to the original optimization problem. This can be achieved iteratively if the Fisher Identity for the surrogate function and the objective function is satisfied. That is, 
\begin{equation}
\left. \frac{\partial}{\partial \bs \theta}\m{J}(\bs \theta) \right \vert_{\bs \theta = \hat{\bs \theta}^{(i)}} = \left. \frac{\partial}{\partial \bs \theta}\m{Q}(\bs \theta,\hat{\bs \theta}^{(i)} ) \right \vert_{\bs \theta = \hat{\bs \theta}^{(i)}}.
\label{eq:Fisher}
\end{equation}
\begin{lemma}
For the class of objective functions in \eqref{eq:cost_kernel}, the surrogate function $\m Q(\bs \theta, \hat{\bs \theta}^{(i)})$ in \eqref{eq:Q_gen} satisfies the Fisher identity defined in \eqref{eq:Fisher}.
\label{lemma:Fisher}
\end{lemma}
\begin{proof}
From \eqref{eq:J_aux} we have: \vspace{-1mm} \small{
\begin{equation*}
 \left. {\frac{\partial}{\partial \bs \theta}\m{J}(\bs \theta)} \right \vert_{\bs \theta \!=\! \hat{\bs \theta}^{(i)}} \!\!=\!\! \left. {\frac{\partial}{\partial \bs \theta} \m{Q}(\bs \theta, \hat{\bs{\theta}}^{(i)})} \right \vert_{\bs \theta = \hat{\bs \theta}^{(i)}} - \left. {\frac{\partial}{\partial \bs \theta} \m{H}(\bs \theta, \hat{\bs{\theta}}^{(i)})} \right \vert_{\bs \theta = \hat{\bs \theta}^{(i)}}.
\end{equation*}}\normalsize
Next, let us consider the gradient of the auxiliary function $\m{H}(\bs \theta,\hat{\bs \theta}^{(i)} )$:
\small{ 
\begin{flalign}
& \left. \frac{\partial}{\partial \bs \theta}\m{H}(\bs \theta,\hat{\bs \theta}^{(i)} ) \right \vert_{\bs \theta = \hat{\bs \theta}^{(i)}} && \nonumber \\ 
&= \bigintsss \limits_{ \Omega(\bs{y})}\! \!\!\!\!\left[ \frac{K(\bs z, \bs \theta)}{\m{V}(\bs \theta)} \right ]^{-1} _{\bs \theta = \hat{\bs \theta}^{(i)}} \!\! \frac{\partial}{\partial \bs \theta} \!\!\left[ \frac{K(\bs z, {\bs \theta})}{\m V({\bs \theta}}\right] _{\bs \theta = \hat{\bs \theta}^{(i)}}\!\! \frac{K(\bs z, \hat{\bs \theta}^{(i)})}{\m V(\hat{\bs \theta}^{(i)})} d\mu(\bs z) && \nonumber \\
& = \!\!\!\bigintsss \limits_{ \Omega(\bs{y})} \frac{\partial}{\partial \bs \theta} \left[ \frac{K(\bs z, {\bs \theta})}{\m V({\bs \theta})}\right] _{\bs \theta = \hat{\bs \theta}^{(i)}}d\mu(\bs z) && \nonumber  \\
&= \frac{\partial}{\partial \bs \theta} \left [  \bigintsss \limits_{ \Omega(\bs{y})}   \frac{K(\bs z, {\bs \theta})}{\m V({\bs \theta})} d\mu(\bs z) \right ] _{\bs \theta = \hat{\bs \theta}^{(i)}} && \nonumber \\ &= 0. &&\nonumber 
\end{flalign}}
\normalsize
Hence, \eqref{eq:Fisher} holds. 
\end{proof}

\begin{remark}
Note that the Fisher identity in Lemma \ref{lemma:Fisher} is well known in the EM-framework. However, we have specialized this result for the problem in this paper (i.e. when $K(\bs z, \bs \theta)$ is not necessarily a probability density function)\eor
\end{remark}
\begin{lemma}
The surrogate function $\m{Q}(\bs \theta,\hat{\bs \theta}^{(i)})$ in \eqref{eq:Q_gen} can be utilized to obtain an adequate surrogate function that satisfies the properties in Mark's approach in \eqref{eq:Marks_tight}-\eqref{eq:Marks_Fisher}.
\end{lemma}
\begin{proof}
Notice that in Mark's approach the optimization problem corresponds to the minimization of the objective function. Hence, to maximize, we have $\m{J}(\bs \theta) = -{\m{Q}}(\bs \theta,\hat{\bs \theta}^{(i)}) + {\m{H}}(\bs \theta,\hat{\bs \theta}^{(i)}) $. From \eqref{eq:J_aux} we can construct the surrogate functions $\tilde{\m{Q}}(\bs \theta,\hat{\bs \theta}^{(i)})$ and $\tilde{\m{H}}(\bs \theta,\hat{\bs \theta}^{(i)})$ since 
\begin{align}
\m{J}(\bs \theta)=-\underbrace{({\m{Q}}(\bs \theta,\hat{\bs \theta}^{(i)}) - {\m{Q}}(\hat{\bs \theta}^{(i)},\hat{\bs \theta}^{(i)}) + \m{J}(\hat{\bs \theta}^{(i)})}_{\tilde{\m Q}(\bs \theta,\hat{\bs \theta}^{(i)})} + \underbrace{({\m{H}}(\bs \theta,\hat{\bs \theta}^{(i)}) - {\m{Q}}(\hat{\bs \theta}^{(i)},\hat{\bs \theta}^{(i)}) + \m{J}(\hat{\bs \theta}^{(i)})}_{\tilde{\m H}(\bs \theta,\hat{\bs \theta}^{(i)})}. 
\label{eq:Q_tilde}
\end{align}
The function $\tilde{\m H}(\bs \theta,\hat{\bs \theta}^{(i)})$ satisfies $\tilde{\m H}(\bs \theta,\hat{\bs \theta}^{(i)})-\tilde{\m H}(\hat{\bs \theta}^{(i)},\hat{\bs \theta}^{(i)})$ $\geq 0$, which implies that $\m J(\bs \theta) \leq \tilde{\m Q}(\bs \theta,\hat{\bs \theta}^{(i)})$, satisfying \eqref{eq:Marks_tight}. From $\tilde{\m Q}(\bs \theta,\hat{\bs \theta}^{(i)}) = {\m{Q}}(\bs \theta,\hat{\bs \theta}^{(i)}) - {\m{Q}}(\hat{\bs \theta}^{(i)},\hat{\bs \theta}^{(i)}) + \m{J}(\hat{\bs \theta}^{(i)})$ we can obtain $\tilde{\m Q}(\hat{\bs \theta}^{(i)},\hat{\bs \theta}^{(i)}) =  \m{J}(\hat{\bs \theta}^{(i)})$, satisfying \eqref{eq:Marks_decr}. Finally, given that the auxiliary function $\tilde{\m H}(\bs \theta,\hat{\bs \theta}^{(i)})$ satisfies \eqref{eq:Fisher}, $\tilde{\m Q}(\bs \theta,\hat{\bs \theta}^{(i)})$ satisfies \eqref{eq:Marks_Fisher}.
\end{proof}
\begin{remark}
Since $\frac{d}{d\bs \theta}\m Q(\bs \theta,\hat{\bs \theta}^{(i)}) = \frac{d}{d\bs \theta} \tilde{\m Q}(\bs \theta,\hat{\bs \theta}^{(i)})$, it is simpler to consider the function $\m Q(\bs \theta,\hat{\bs \theta}^{(i)})$ instead of $\tilde{\m Q}(\bs \theta,\hat{\bs \theta}^{(i)})$ in penalized (regularized) and MAP estimation problems, as shown in \textup{\cite{ref:Carvajal2015}} and \textup{\cite{Godoy2014}}.\eor
\end{remark}
We summarize our proposed algorithm in Table \ref{table:Marks}. 
\begin{table}[t]
\caption{Proposed algorithm.} \vspace{0.5mm}  
\centering                          
\setlength{\tabcolsep}{2.5pt}
\begin{tabular}{r l}             
\toprule
Step 1: & Find a kernel that satisfies \eqref{eq:cost_kernel}.\\
{Step 2:} & $i=0$.\\
{Step 3:} &  Obtain an initial guess $\hat{\bs{\theta}}^{(i)}$. \\ 
{Step 4:} &  Compute $\m{Q}(\bs \theta, \hat{\bs \theta}^{(i)})$ as in \eqref{eq:Q_gen}. \\ 
{Step 5:} &  Compute $\tilde{\m{Q}}(\bs \theta, \hat{\bs \theta}^{(i)})$.\\
{Step 6:} &  Incorporate $\tilde{\m{Q}}(\bs \theta, \hat{\bs \theta}^{(i)})$ in the optimization problem and solve.\\
{Step 7:} &  $i = i+1$ and back to Step 4 until convergence.\\
\bottomrule
\end{tabular}

\label{table:Marks}
\end{table}

\section{A quadratic surrogate function for a class of kernels}\label{section:quadratic}

In this paper we focus on a special class of the kernel functions $K(\bs z, \bs \theta)$. For this particular class, the following is satisfied:
\begin{equation}
\frac{\partial}{\partial \bs \theta}\log [K(\bs z, \bs \theta)] = \textbf{A}(\bs z)\bs \theta + \textbf{b}, 
\label{eq:quad_kernel}
\end{equation}
where $\textbf{A}(\bs z)$ is a matrix and $\textbf{b}$ is a vector, both of adequate dimensions. Then, we have that
\begin{align}
& \frac{\partial}{\partial \bs \theta}{\m Q}(\bs \theta, \hat{\bs \theta }^{(i)} ) = \bigintsss \limits_{ \Omega(\bs{y})} [\textbf{A}(\bs z)\bs \theta + \textbf{b}] \frac{K(\bs z, \hat{\bs \theta}^{(i)})}{\m V(\hat{\bs \theta}^{(i)})} d\mu(\bs z) \nonumber \\
& = \left[ \bigintsss \limits_{ \Omega(\bs{y})} \!\!\!\!\! \textbf{A}(\bs z)\frac{K(\bs z,\hat{\bs \theta}^{(i)})}{\m{V}(\hat{\bs \theta}^{(i)})}d\mu(\bs z)  \right]\!\! \bs{\theta} \nonumber + \textbf{b} \!\!\bigintsss \limits_{ \Omega(\bs{y})}\!\! \frac{K(\bs z,\hat{\bs \theta}^{(i)})}{\m{V}(\hat{\bs \theta}^{(i)})}d\mu(\bs z)\nonumber \\
& = \textbf{R}\bs{\theta} + \textbf{b}.
\end{align}
\begin{remark}
Notice that the previous expression is linear with respect to $\bs \theta$. This implies that the function ${\m Q}(\bs \theta, \hat{\bs \theta }^{(i)} )$  is quadratic with respect to the parameter vector $\bs \theta$. \eor
\end{remark}
From the Fisher Identity in \eqref{eq:Fisher} we have that
\begin{equation}
\left. \frac{\partial}{\partial \bs \theta}\m{J}(\bs \theta) \right \vert_{\bs \theta = \hat{\bs \theta}^{(i)}} = \textbf{R}\hat{\bs \theta}^{(i)} + \textbf{b},
\end{equation}
from which we can solve for \textbf{R} in some cases\footnote{The matrix \textbf{R} can also be computed using Monte Carlo algorithms.}. In particular, if $\textbf{A}(\bs z)$ is a diagonal matrix, then \textbf{R} is also a diagonal matrix defined by $\textbf{R} = \text{diag}[r_1, r_2, ...]$. Thus, we have
\[
\left. \frac{\partial}{\partial \theta_k}\m{J}(\bs \theta) \right \vert_{\bs \theta = \hat{\bs \theta}^{(i)}} \!\!\! = r_k \hat{\theta}^{(i)}_i \!+ b_k  \Rightarrow r_k \!=\! \frac{\left. \frac{\partial}{\partial \theta_k}\m{J}(\bs \theta) \right \vert_{\bs \theta = \hat{\bs \theta}^{(i)}} - b_k}{\hat{\theta}_k^{(i)}}, 
\]
where $\theta_i$ is the $i$th component of the parameter vector $\bs \theta$, $\hat{\bs \theta}_i^{(i)}$ is the $i$th component of the estimate $\hat{\bs \theta}^{(i)}$, $r_i$ is the $i$th element of the diagonal of $\textbf{R}$, and $b_i$ is the $i$th element of the vector $\textbf{b}$.
Hence, when optimizing the auxiliary function $\m{Q}(\bs  \theta, \hat{\bs \theta}^{(i)})$ we obtain
\[ \frac{\partial}{\partial \bs \theta}\m{Q}(\bs  \theta, \hat{\bs \theta}^{(i)}) = \begin{bmatrix}
r_1 & & \\ & r_2 & \\ & & \ddots
\end{bmatrix} \bs \theta + \textbf{b} = 0 \, \Rightarrow \, \hat{\theta}_k^{(i+1)} = - \frac{b_k}{r_k}.\]
Equivalently, 
\begin{equation}
\hat{\bs \theta}^{(i+1)} = \textbf{R}^{-1}\textbf{b}.
\end{equation}
This implies that in our approach, it is not necessary to obtain the auxiliary function $\m Q$ and optimize it. Instead, by computing $\textbf{R}$ and $\textbf{b}$ at every iteration, the new estimate can be obtained.

For the class of kernels here described, the proposed method for constructing surrogate functions can also be understood as part of sequential quadratic programming (SQP) methods \cite[ch. 12.4]{Fletcher1987}. Indeed, the general case of equality and inequality-constrained minimization problems is defined as \cite[ch. 4]{Izmailov2014}:
\begin{align*}
\bs \theta^{*} = & \, \text{arg} \min_{\bs \theta}  f(\bs \theta) \\
\text{s. t.   }  & \,\, h(\bs \theta) = 0 \, , \, g(\bs{\theta}) \leq 0,
\end{align*}
which is solved by iteratively defining quadratic functions that approximate the objective function and the inequality constraint around a current iterate $\hat{\bs \theta}^{(i)}$. In the same way, our proposal generates an algorithm with quadratic surrogate functions.

\section{Numerical Examples}

\label{section:example}
In this section, we illustrate our proposed algorithm with two numerical examples.
\subsection{Example 1: Estimation of the Distribution of Stellar Rotational Velocities}

One of the many problems in Astronomy deals with is the estimation of rotational velocities of stars. This particular problem is of great importance, since it allows astronomers to describe and model the stars formation, their internal structure and evolution, as well as how they interact with other stars, see e.g. \cite{ref:Christen2016, ref:Cure2014, ref:Chandrasekhar1950}. 

Modern telescopes allow for the measurement of the rotational velocities from the telescope point of view, that is, a projection of the true rotational velocity. This is modelled (spatially) as the convolution of the true rotational velocity pdf and a uniform distribution over the sphere (for more details see e.g. \cite{ref:Cure2014}):
\begin{equation}
p(y|\sigma) = \bigintssss_y^\infty \frac{y}{x\sqrt{x^2 - y^2}}p(x|\sigma)dx,
\label{eq:conv_star}
\end{equation}
where $p(y|\sigma)$ is the uniform projected rotational velocity pdf and $p(x|\sigma)$ is the true rotational velocity pdf to be estimated, and $\sigma$ a hyperparameter.

A common model for $p(x|\sigma)$ found in the Astronomy literature is the Maxwellian distribution (see e.g. \cite{ref:Cure2014, ref:Deutsch1969})
\begin{equation}
p(x|\sigma) = \sqrt{\frac{2}{\pi}}\frac{1}{\sigma^3}x^2e^{-{x^2}/(2\sigma^2)}.
\end{equation}
In practice, the measurements correspond to realizations of $p(y|\sigma)$ \cite{ref:Cure2014}, from which the likelihood function can be defined as:
\begin{equation}
p(\textbf{y}|\sigma) = \prod_{k=1}^N p(y_k|\sigma),
\end{equation}
where $\textbf{y} = [y_1, ... , y_N]^T$,
\[p(y_k|\sigma) = \bigintssss_{y_k}^\infty \frac{y_k}{x_k\sqrt{x_k^2 - y_k^2}}p(x_k|\sigma)dx_k,\]
$x_k$ is Maxwellian distributed, and $N$ is the number of measurement points. Hence, the log-likelihood function can be expressed as:
\begin{equation}
\ell(\sigma) =\sum_{k=1}^{N} \log  \left[\bigintssss_{y_k}^{\infty} \!\! \frac{y_k}{x_k\sqrt{x_{k}^2-y_k^2}}p(x_k|\sigma)dx_k\right].
\label{eq:log_lik_star}
\end{equation}
If we define the complete data $\bs z = (\textbf{x},\textbf{y})$, the kernel function $K(\cdot,\cdot)$ and the measure $\mu(\cdot)$ in \eqref{eq:cost_kernel} can be defined as
\begin{equation}
K(x_k,\sigma)=p(x_k|\sigma)=\sqrt{\frac{2}{\pi}}\frac{x_k^2}{\sigma^3}e^{-{x_k^2}/(2\sigma^2)},
\label{eq:K_Maxwellian}
\end{equation}
and
\begin{equation}
d\mu(x_k,y_k)=\frac{y_k}{x_k\sqrt{x_k^2-y_k^2}}dx_k.
\end{equation}
Then, the log-likelihood function in \eqref{eq:log_lik_star} can be written as:
\begin{equation}
\ell(\sigma)=\sum_{k=1}^{N}\log\left[\mathcal{V}_k(\sigma)\right],
\end{equation} 
with
\begin{equation}
\mathcal{V}_k(\sigma)=\bigintssss_{y_k}^{\infty}K(x_k,\sigma)d\mu(x_k,y_k),
\end{equation}
Thus, the ML estimator is obtained from:
\begin{equation}
\hat{\sigma}_\text{ML}=\arg \max_{\sigma}\; \sum_{k=1}^{N}\log\mathcal{V}_k(\sigma).
\label{eq:MLE_star}
\end{equation}
Since the parameter that is needed to be estimated is part of the convolution in \eqref{eq:conv_star}, the optimization problem in \eqref{eq:MLE_star} cannot be solved in a straightforward manner. Instead, we utilize the re-interpretation of the EM algorithm that we propose for solving \eqref{eq:MLE_star}.

First, notice that from the surrogate function $\bar{\mathcal{Q}}(\sigma,\hat{\sigma}^{(i)})$ can be expressed as:
\begin{equation}
{\mathcal{Q}}(\sigma,\hat{\sigma}^{(i)})=\sum_{k=1}^{N} \mathcal{Q}_k(\sigma,\hat{\sigma}^{(i)}),
\label{eq:Q_star}
\end{equation}
with
\begin{equation}
\mathcal{Q}_k(\sigma,\hat{\sigma}^{(i)}) \! = \!\! \bigintssss_{y_k}^{\infty}\! \! \!\! \! \log (K(x_k,\sigma))\frac{K(x_k,\hat{\sigma}^{(i)})}{\mathcal{V}_k(\hat{\sigma}^{(i)})}d\mu(x_k,y_k).
 \end{equation}
For convenience, we can differentiate the auxiliary function $\m{Q}(\sigma,\hat{\sigma}^{(i)})$ in \eqref{eq:Q_star} with respect to $1/\sigma$ obtaining:
\begin{equation}
\frac{\partial \bar{\mathcal{Q}}(\sigma,\hat{\sigma}^{(i)})}{\partial(1/\sigma)}=\sum_{k=1}^{N}\bigintssss_{y_k}^{\infty}\!\left[3\sigma - \frac{x_k^2}{\sigma}\right]\!\!\frac{K(x_k,\hat{\sigma}^{(i)})}{\mathcal{V}_k(\hat{\sigma}^{(i)})}d\mu(x_k,y_k). 
\end{equation}
Then, equating to zero and solving for $\sigma$ we finally obtain
\begin{equation}
\hat{\sigma}^{(i+1)} = \sqrt{\frac{\mathcal{S}(\textbf{y},\hat{\sigma}^{(i)})}{3N}},
\label{eq:sigma_update}
\end{equation}
where
\begin{equation}
\label{eq:S(xt,beta)}
\m{S}(\textbf{y},\hat{\sigma}^{(i)})=\sum_{t=1}^{N}\bigintssss_{y_k}^{\infty}x_k^2\frac{K(x_k,\hat{\sigma}^{(i)})}{\mathcal{V}_k(\hat{\sigma}^{(i)})}d\mu(x_k,y_k).
 \end{equation}

In Table \ref{Table:algorithm} we summarized the specialisation of our proposed algorithm for this example.
 \begin{table}[h!]
 	\caption{Proposed algorithm for Maxwellian distribution estimation in Example 1.}
 	\centering                          
\setlength{\tabcolsep}{2.5pt}
 	\label{Table:algorithm}
 	\begin{tabular}{cl}
 		\toprule
 		{Step 1:} & $i=0$.\\
{Step 2:} &  Obtain an initial guess $\hat{\bs{\sigma}}^{(i)}$. \\ 
 		Step 3: & Compute the integral given by (\ref{eq:S(xt,beta)}).  \\
 		Step 4: & Compute $\hat{\sigma}^{(i+1)}$ \eqref{eq:sigma_update}\\
 		Step 5: & $i=i+1$ and back to Step 3 until convergence.\\ \bottomrule
 	\end{tabular}
 \end{table}

\begin{figure}[h!]
\centering
\includegraphics[width=0.55\textwidth, trim=10mm 2mm 17mm 9mm, clip]{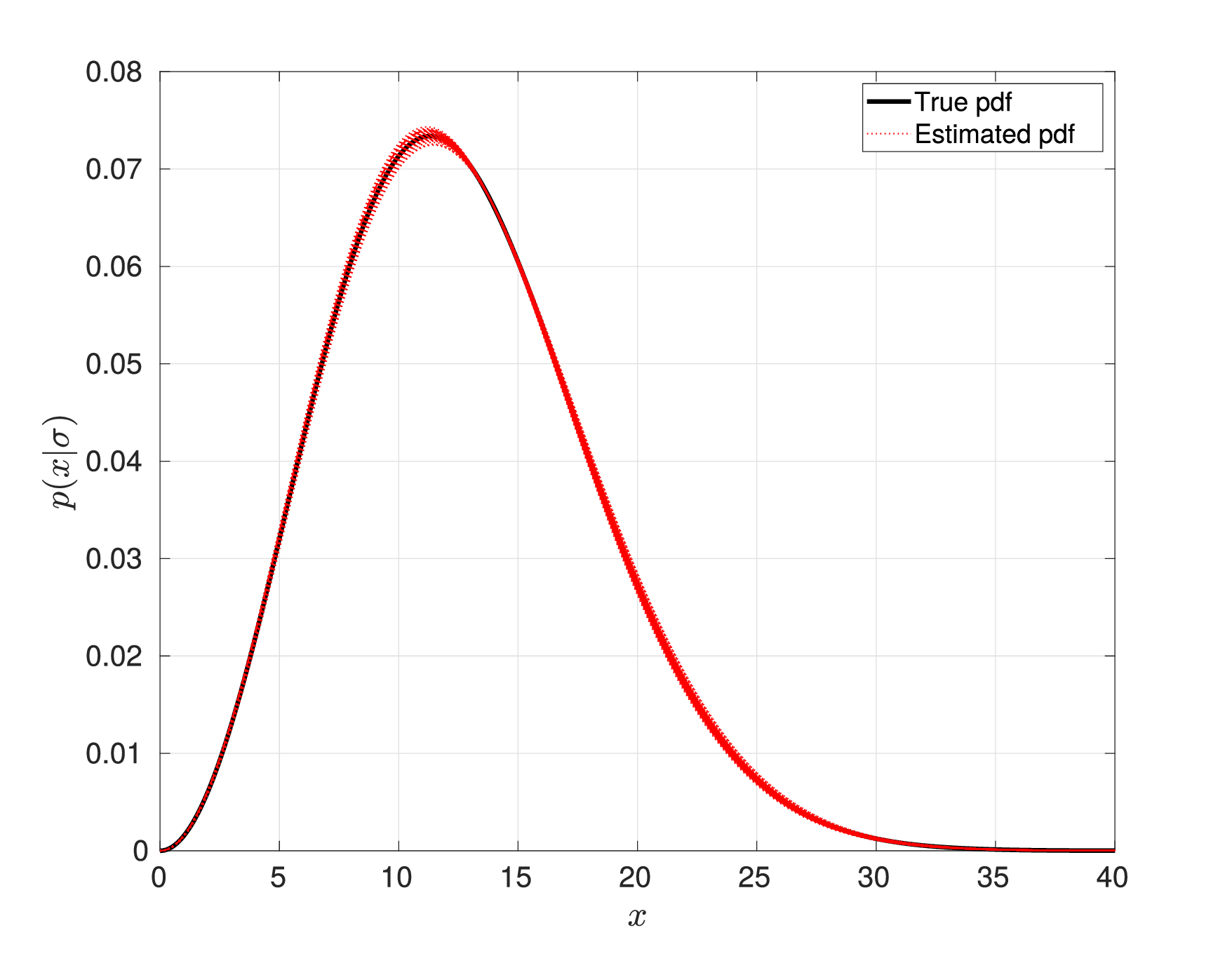}
\caption{Convergence of the proposed approach to the global optimum}
\label{fig:example1_1}
\end{figure}

\begin{figure}[h!]
\centering
\includegraphics[width=0.55\textwidth, trim=10mm 2mm 17mm 12mm, clip]{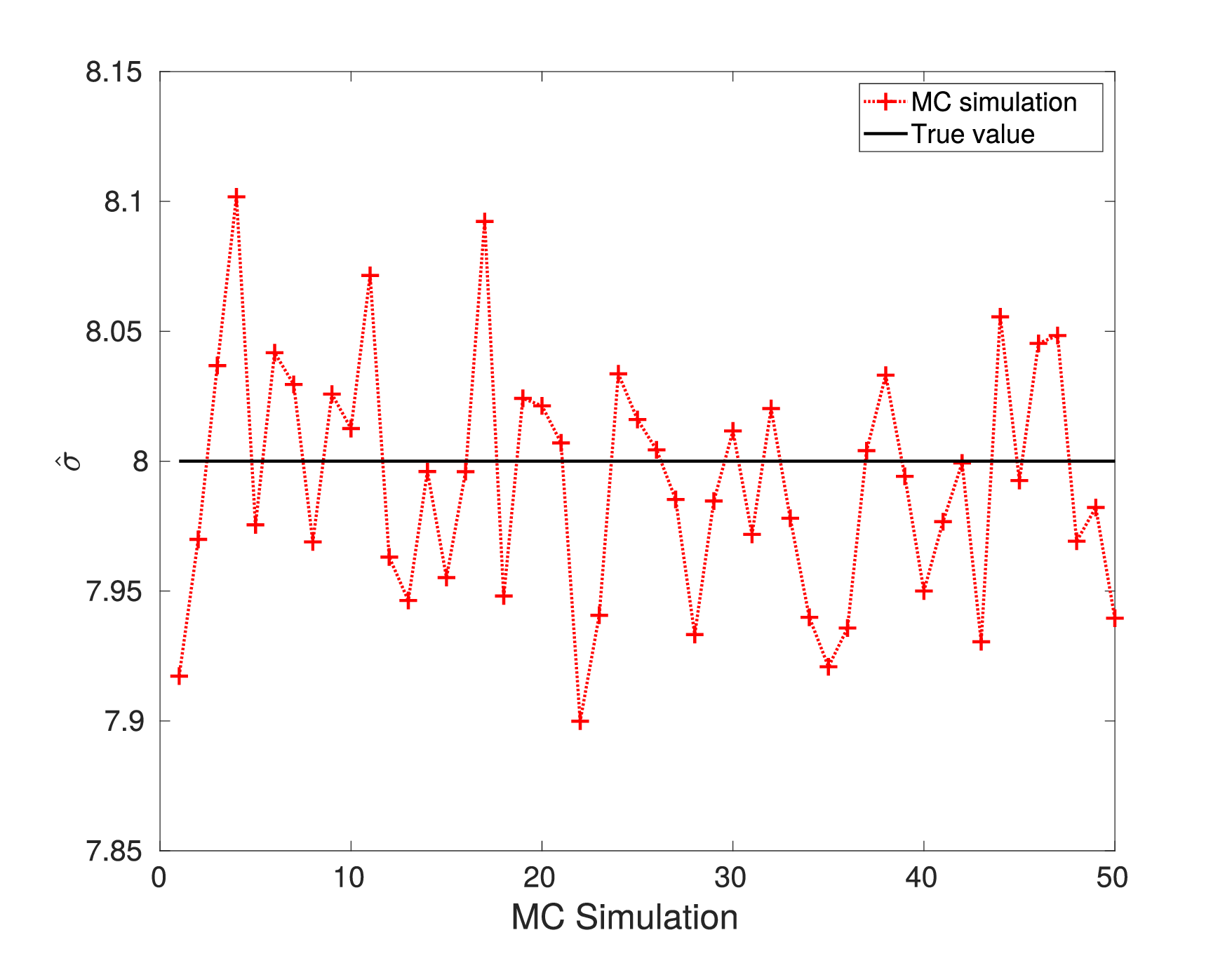}
\caption{Convergence of the proposed approach to the global optimum}
\label{fig:example1_2}
\end{figure}

For the numerical simulation, we have considered the problem solved in \cite{ref:Cure2014}, with the true dispersion parameter $\sigma_0 = 8$. The measurement data $\textbf{y} = [y_1,...,y_N]$ was generated using the \textit{Slice Sampler} (see e.g. \cite{ref:Robert}) applied to \eqref{eq:conv_star}. The simulation setup is as follows:
\begin{itemize}
	\item The data length is given by $N=10000$.
	\item The number of Monte Carlo (MC) simulations is 50.
	\item The stopping criterion is given by:
	\begin{equation*}
	\Vert \hat{\sigma}^{(i)}-\hat{\sigma}^{(i-1)}\Vert / \Vert \hat{\sigma}^{(i)}\Vert < 10^{-6},
	\end{equation*}
	or the maximum number of iterations of $100$ has been reached.
\end{itemize}
The results are shown in Fig. \ref{fig:example1_1}, were the estimated $p(x|\sigma)$ for each MC simulation is shown. It is clear that the estimated Maxwellian distributions are very similar to the \textit{true} density distribution. The mean value of the estimated parameter was $\hat{\sigma}=7.9920$. The estimation from each MC simulation is shown in Fig \ref{fig:example1_2}. It can be clearly seen that the estimated parameter $\hat{\sigma}$ is close to the \textit{true} value.

 \vspace{-3mm}

\subsection{Example 2: Inequality Constrained ML Estimation}

Let us consider the following inequality constrained ML estimation problem:
\begin{align}
\bs \theta^{*} = & \, \text{arg} \min_{\bs \theta}  \Vert \bs{y} - \bs{A \theta} \Vert_2^2  \nonumber \\
\text{s. t.   }  & f(\bs \theta) = \sum_{m = 1}^M \left \Vert \frac{\sqrt{\beta_m}}{\tau} \bs{\theta}_m \right \Vert_2^q \leq \gamma,
\label{eq:example}
\end{align}
where $\bs \theta = [\bs \theta_1^T, ..., \bs \theta_M^T]^T$ is an unknown vector parameter, $\bs A$ is the regressor matrix, and $\bs y$ corresponds to the measurements, with $q = 0.4$, $\tau = 0.2$,  and $\gamma = 7$. Notice that the constraint function belongs to a family of sparsity constraints found in sparse system identification and optimization -- see e.g \cite{ref:Hyder2009} and the references therein. The function of the unknown parameter in the inequality constraint can be understood as its ``group $\ell_q$-norm'', with $M$ groups and where $\bs \theta_m$ is the $m$th group of length $\beta_m$\footnote{Notice that when $m=1$ we obtain the standard $\ell_q$-norm, and when $q = 1$ we obtain the $\ell_1$-norm utilized in the Lasso.}. Hence, the inequality constraint can be expressed as a Multivariate Power Exponential (MPE) distribution \cite{ref:Solaro2004} of $\bs \theta$. On the other hand, the MPE distributions can be expressed as VMGMs \cite{ref:Lange1993}, which allows for utilizing the ideas presented in Section \ref{section:Gen_app}, where the corresponding kernel function is a multivariate Gaussian distribution. In this example, the elements of the matrix $\bs A$ were generated using a Normal distribution ($\m N(0,1)$) and the measurement data was generated from
\[ \bs y = \bs{A \theta} + \bs n,  \]
where $\bs n \sim \m N(\bs 0, 0.01\textbf{I})$ is the additive measurement noise. Here the number of measurements is $256$ and $\bs \theta = [0.54,$ $1.83, -2.26,$ $0,0, 0, 0, 0,0,0.86,0.32,-1.31]^T$, with $\beta_m = 3$ $\forall m$, and thus $M = 4$. We have also considered $150$ Monte Carlo simulations. 

For the attainment of the surrogate function we consider the following:
\begin{enumerate}
\item When expressed as a MPE distribution, the group $\ell_q$-norm is given by \cite{ref:Solaro2004},
\[ \hspace{-5mm}p(\bs \theta_m) = \kappa_m(q,\tau) e^{-\left \Vert \frac{\sqrt{\beta_m}}{\tau} \bs{\theta}_m \right \Vert^q_2}, \]
where $\kappa_m(q,\tau)$ is a constant that depends on $\tau$ and $q$, and 
\[p(\bs \theta) = \prod_{m = 1}^M p(\bs \theta_m).\]
\item When expressed as a VMGM, the group $\ell_q$-norm is given by
\[ p(\bs \theta_m) = \int \m{N}_{\bs \theta}(0,\bs \Sigma(\lambda_m))p(\lambda_m)d\lambda_m. \]
\item With the previous expressions, and using the results in Section \ref{section:quadratic}, we notice that
\[\frac{\partial \m Q(\bs \theta, \hat{\bs \theta}^{(i)})}{\bs \theta_m} = -k_m^{(i)} \bs \theta_m,\]
where $\hat{\bs \theta}_m^{(i)}$ is the estimate of the $m$th group at the $i$th iteration and $k_m^{(i)}$ is a constant that depends on $\hat{\bs \theta}_m^{(i)}$.
\item Then we obtain $\m Q(\bs \theta, \hat{\bs \theta}^{(i)}) = \frac{k_m^{(i)}}{2} \left( \bs \theta_m^T \bs \theta_m \right) + c$, where $c$ is a constant.
\end{enumerate}

\begin{figure}[t!]
\centering
\includegraphics[width=0.65\textwidth, trim=43mm 2mm 22mm 8mm, clip]{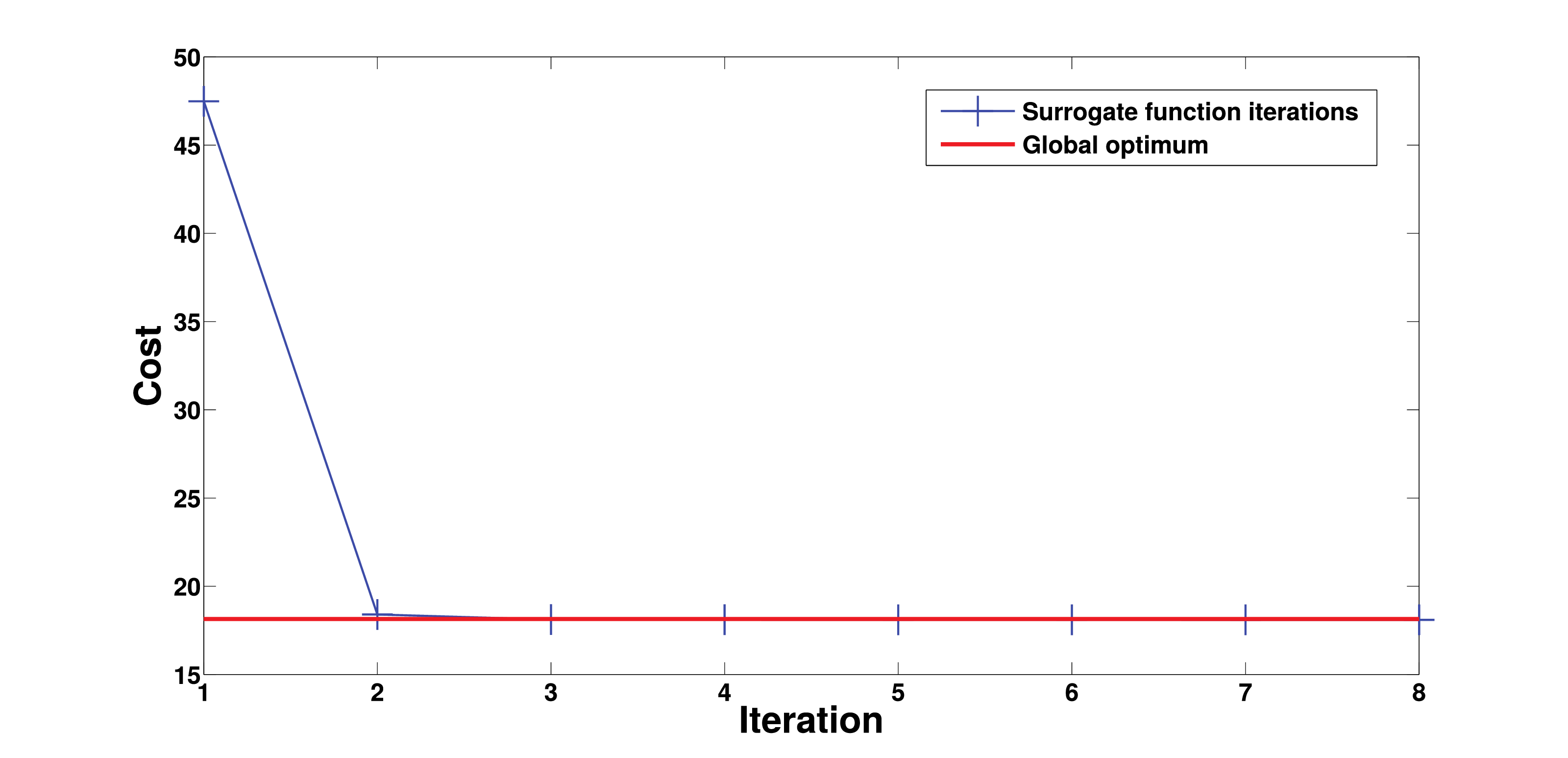}
\caption{Convergence of the proposed approach to the global optimum}
\label{fig:example2}
\end{figure}

Then, the inequality constraint in \eqref{eq:example} can be replaced with the following constraint (that includes the surrogate function $\tilde{\m Q}(\bs \theta, \hat{\bs \theta}^{(i)})$ as shown in \eqref{eq:Q_tilde})
\begin{equation}
\sum_{m=1}^M  \frac{k_m^{(i)}}{2} \left( \bs \theta_m^T \bs \theta_m - (\hat{\bs \theta}_m^{(i)})^T (\hat {\bs \theta}_m^{(i)})  \right) + \left \Vert \frac{\sqrt{\beta_m}}{\tau} \hat{\bs{\theta}}_m^{(i)} \right \Vert_2^q \leq \gamma,
\end{equation}
The results of the optimization problem for one realization is shown in Fig. \ref{fig:example2}, where the optimization problem with the surrogate function was solved using the optimization software CVX \cite{ref:cvx1, ref:cvx2} and the original non-convex problem was solved using the optimization software BARON \cite{ref:baron1, ref:baron2}, which allowed us to obtain the global optimum. From Fig. \ref{fig:example2} it can be seen that our approach finds the global minimum with a few iterations only. 
On the other hand, we have also considered the mean square error (MSE) of the estimates and compared them with the unconstrained problem (least squares). The MSE obtained with our approach is $3.6\times 10^{-2}$, whilst the least squares estimate yielded an MSE equal to $4.34\times10^{-2}$. Clearly, the incorporation of the constraint aided the estimation.

\section*{Conclusions}
\label{section:conclusions}
In this paper we have presented a systematic approach for constructing surrogate functions in a wide range of optimization problems. Our approach can be utilized for constructing surrogate functions for both the cost function and the constraints, generalizing the popular EM and MM algorithms. Our approach is based on the utilization of data augmentation and kernel functions, yielding simple optimization algorithms when the kernel can be expressed as VMGM. 

\section*{Acknowledgments}
This work was partially supported by FONDECYT~-~Chile through grants No. 3140054 and 1181158.  This work was also partially supported by the Advanced Center for Electrical and Electronic Engineering (AC3E, Proyecto Basal FB0008), Chile. The worl of R. Orellana was partially supported by PIIC program, scholarship 015/2018, UTFSM.

\bibliographystyle{plain}   
\bibliography{biblio}

\end{document}